\font\we=cmb10 at 14.4truept
\font\li=cmb10 at 12truept
\noindent
\centerline {\we A Note on Arithmetic Cohomologies for Number Fields}
\vskip 0.45cm
\centerline {\li Lin WENG}
\vskip 0.30cm
\centerline {\bf Graduate School of Mathematics, Nagoya University, Japan}
\vskip 0.50cm
\noindent
In this short article, we develop a new cohomology theory for
number fields as a part of our Program for Geometric Arithmetic.
\vskip 0.30cm
\noindent
{\li 1. Arithmetic Cohomology Groups}
\vskip 0.30cm
Let $F$ be a number field with discriminant $\Delta_F$. Denote its
(normalized) absolute values by
$S_F$, and write $S_F=S_{\rm fin}\buildrel \cdot\over\cup S_\infty$,
where $S_\infty$ denotes the collection of all archimedean valuations.
For simplicity, we use $v$ (resp. $\sigma$) to denote elements in $S_{\rm
fin}$ (resp. $S_\infty$).

Denote by ${\bf A}={\bf A}_F$ the ring of adeles of $F$, and ${\rm
Gl}_r({\bf A})$ the rank $r$ general linear group over {\bf A}, and write
${\bf A}:={\bf A}_{\rm fin}\oplus {\bf A}_\infty$ and ${\rm Gl}_r({\bf
A}):={\rm Gl}_r({\bf A})_{\rm fin}\times {\rm Gl}_r({\bf A})_\infty$
according to their finite and infinite parts.

For any $g=(g_{\rm fin};g_\infty)=(g_v;g_\sigma)\in {\rm Gl}_r({\bf A})$,
define the injective morphism $i(g):=i(g_\infty):F^r\to {\bf
A}_\infty^r$ by $(f)\mapsto (g_\sigma\cdot f)$, $F^r(g):={\rm
Im}\big(i(g)\big)$ and
$${\bf A}^r(g):=\{(a_v;a_\sigma)\in {\bf A}^r:\exists f\in F^r\ {\rm
s.t.}\ g_v^{-1}(a_v)\in {\cal O}_v^r,\ g_v^{-1}(f)\in {\cal
O}_v^r,\,\forall v\ {\rm and}\ (a_\sigma)=i(g_\infty)(f)\}.$$ Then we
have the following 9-diagram with all columns and roots exact:
$$\matrix{&&0&&0&&0&&\cr
&&\downarrow&&\downarrow&&\downarrow&&\cr
0&\to&{\bf A}^r(g)\cap F^r(g)&\to&{\bf A}^r(g)&\to&{\bf A}^r(g)/{\bf
A}^r(g)\cap F^r(g)&\to&0\cr
&&\downarrow&&\downarrow&&\downarrow&&\cr
0&\to& F^r(g)&\to&{\bf A}^r&\to&{\bf A}^r/F^r(g)&\to&0\cr
&&\downarrow&&\downarrow&&\downarrow&&\cr
0&\to& F^r(g)/{\bf A}^r(g)\cap F^r(g)&\to&{\bf A}^r/{\bf A}^r(g)&\to&{\bf
A}^r/{\bf A}^r(g)+F^r(g)&\to&0\cr
&&\downarrow&&\downarrow&&\downarrow&&\cr
&&0&&0&&0.&&\cr}$$

Motivated by this and Weil's adelic cohomology theory for divisors over
algebraic curves, (see e.g., [Serre] and [Weil]), we introduce the
following
\vskip 0.30cm
\noindent
{\bf Definition.} {\it For any $g\in {\rm GL}_r({\bf A})$, define its
0-th and 1-st arithmetic cohomology groups by
$$H^0(F,g):={\bf A}^r(g)\cap F^r(g),\qquad{\rm and}\qquad H^1(F,g):={\bf
A}^r/{\bf A}^r(g)+F^r(g).$$}

\noindent
{\bf Theorem.} (Serre Duality=Pontrjagin Duality) {\it As locally compact
groups, $$H^1(F,g)\simeq \widehat{H^0(F,k_F\otimes g^{-1})}.$$ Here $k_F$
denotes an idelic canonical element of $F$.}
\vskip 0.30cm
\noindent
{\it Remark.} For many $v\in S_{\rm fin}$, denote $\partial_v$ the local
differential of $f_v$, the $v$-completion of $F$ at $v$. Denote by
${\cal O}_v$ the valuation ring and $\pi_v$ any local parameter. Then
$\partial_v=\pi_v^{{\rm ord}_v(\partial_v)}\cdot {\cal O}_v$ and we call
$k_F:=(\partial_v^{-{\rm ord}_v(\partial_v)};1)$ an idelic canonical
element of $F$.
\vskip 0.30cm
\noindent
{\it Proof.} As usual, introduce a basic character $\chi$ on {\bf A} by
$\chi=(\chi^{(r)}_v;\chi^{(r)}_\sigma)$ where $\chi_v:=\lambda_v\circ {\rm
Tr}^{F_v}_{{\bf Q}_v}$ with $\lambda_v:{\bf Q}_v\to {\bf Q}_v/{\bf
Z}_v\hookrightarrow {\bf Q}/{\bf Z}\hookrightarrow {\bf R}/{\bf Z}$, and
$\chi_\sigma:=\lambda_\infty\circ {\rm Tr}_{\bf R}^{F_\sigma}$ with
$\lambda_\infty:{\bf R}\to {\bf R}/{\bf Z}$.  Then the pairing
$(x,y)\mapsto e^{2\pi i\chi(x\cdot y)}$ induces a natural isomorphism
$\widehat {{\bf A}^r}\simeq {\bf A}^r$. Moreover, with respect to this
pairing,
  $\big({\cal O}_v^r\Big)^\perp=\partial_v^r$ and $F^r(g)^\perp\simeq
F^r(g^{-1})=F^r(k_F\otimes g^{-1})$. Thus, we could complete the proof by
the fact that
$$\Big({\bf A}^r(g)\cap F^r(g)\Big)^\perp=\big({\bf A}^r(g)\big)^\perp
+\big(F^r(g)\big)^\perp.$$
\vskip 0.30cm
\noindent
{\li 2. Arithmetic Counts}
\vskip 0.30cm
Motivated by the  fact that the dimension of a vector space is
equal to the dimension of its dual space, one of the basic principal we
adopt in counting locally compact groups is the following:
\vskip 0.30cm
\noindent
{\bf Counting Axiom.} {\it If $\#_{\rm ga}$ is a count for a certain
class of locally compact groups $G$,  then $\#_{\rm ga}(G)=\#_{\rm
ga}(\hat G)$.}
\vskip 0.30cm

Clearly, this is compactible with the Pontrjagin duality. Thus, our
counts of arithmetic cohhomology groups should be based on
Fourier analysis over these groups, or better, the Fourier inverse
formula for Fourier transforms. In this way, practically, it is natural
to use the Plancherel formula to do counts.

While we may use any reasonable test functions on ${\bf A}^r$ to do the
counts, for simplicity and also as a continuation of a more classical
mathematics, we set $f:=\prod_v f_v\cdot\prod_\sigma f_\sigma$.
Here $f_v$ is the characteristic function of ${\cal O}_v^r$;
$f_\sigma(x_\sigma):=e^{-\pi|x_\sigma|^2/2}$ if $\sigma$ is real; and
$f_\sigma(x_\sigma):=e^{-\pi|x_\sigma|^2}$ if $\sigma$ is complex.
\vskip 0.30cm
\noindent
{\bf Definition.} {\it (1) The arithemtic counts of the 0-th and the
1-st arithmetic cohomology groups for  $g\in {\rm
Gl}_r({\bf A})$ are defined to be
$$\eqalign{\#_{\rm ga}\big(H^0(F,g)\big):=&\#_{\rm ga}\Big(H^0(F,g); f,
dx\Big):=\int_{H^0(F,g)}|f(x)|^2dx;\cr
\#_{\rm ga}\big(H^1(F,g)\big):=&\#_{\rm ga}\Big(H^1(F,g); \hat f,
d\xi\Big):=\int_{H^1(F,g)}|\hat f(\xi)|^2d\xi.\cr}$$ Here $dx$ denotes
(the restriction of) the standard Haar measure on {\bf A},  $d\xi$
denotes (the induced quotient measure from) the dual measure, and $\hat
f$ denotes the Fourier transform of $f$;
\vskip 0.30cm
\noindent
(2) The 0-th and the 1-st arithmetic cohomologies of $g\in {\rm
Gl}_r({\bf A})$ is defined to be
$$h^0(F,g):=\log \Big(\#_{\rm ga}\big(H^0(F,g)\big)\Big)\qquad{\rm
and}\qquad h^1(F,g):=\log \Big(\#_{\rm ga}\big(H^1(F,g)\big)\Big).$$}

\noindent
{\li 3. Serre Duality and Riemann-Roch}
\vskip 0.30cm
For the arithmetic cohomologies just introduced, we have the following
\vskip 0.30cm
\noindent
{\bf Theorem.} (1) (Serre Duality) $h^1(F,g)=h^0(F,k_F\otimes g)$;
\vskip 0.30cm
\noindent
(2) (Riemann-Roch Theorem)
$$h^0(F,g)-h^1(F,g)={\rm deg}(g)-{r\over 2}\cdot\log|\Delta_F|.$$

\noindent
{\it Proof.} (1) is a direct consequence of the topological version of
Serre duality and the Plancherel Formula.
\vskip 0.25cm
(2) is a direct consequence of the Serre duality just proved and the
Poisson summation formula by noticing that $\Big(H^0(F,g)\Big)^\perp
=H^0(F,k_F\otimes g^{-1})$. This then completes the proof.
\eject\vskip 0.30cm
\noindent
{\li 4. Comments}
\vskip 0.30cm
Our work here is motivated by the works of Weil, Tate, van der Geer and
Schoof, and Li, as well as the works of Lang, Arakelov, Szpiro, Neukirch,
Connes, Deninger, Borisov, and Moreno. For details, please see the
references below, in particular [Weng2]. As an application, we could
introduce new yet natural non-abelian zeta functions for number fields,
based on a discussion about intersection stability ([Weng1]).
\vskip 0.30cm
As it stands, one may start from only ${\bf A}_\infty$ to introduce
arithmetic cohomogy groups and do corresponding counts. This has the
advantage of being more concrete, and applies well say in Arakelov theory.
On the other hand, a more general theory indeed works well for a much
wider class of characters and  test functions. We leave this to the
reader. (See e.g., [Tate] and [Moreno].)
\vskip 0.30cm
Finally, I would like to thank Deninger for his discussion
and encouragement.
\vskip 1.50cm
\centerline {\li REFERENCES}
\vskip 0.25cm
\item{[Borisov]} A. Borisov, Convolution structures and arithmetic
cohomology, preprint
\vskip 0.25cm
\item{[Connes]} A. Connes, Trace formula in
noncommutative geometry and the zeros of the Riemann zeta function.  Sel.
math. New ser  {\bf 5}  (1999),  no. 1, 29--106.
\vskip 0.25cm
\item{[Deninger]} Ch. Deninger, Motivic $L$-functions and regularized
determinants,  Motives (Seattle, WA, 1991), 707-743, Peoc. Sympos. Pure
Math, {\bf 55} Part 1, AMS, Providence, RI, 1994
\vskip 0.25cm
\item{[GS]} G. van der Geer \& R. Schoof, Effectivity of Arakelov
Divisors and the Theta Divisor of a Number Field, Sel. Math., New ser.
{\bf 6} (2000), 377-398
\vskip 0.25cm
\item{[Lang1]} S. Lang, {\it Algebraic Number Theory},
Springer-Verlag, 1986
\vskip 0.25cm
\item{[Lang2]} S. Lang, {\it Fundamentals on Diophantine Geometry},
Springer-Verlag, 1983
\vskip 0.25cm
\item{[Li]} X. Li, A note on the Riemann-Roch theorem for function fields.
{\it Analytic number theory}, Vol. {\bf 2},  567--570, Progr. Math.,
{\bf 139}, Birkh\"auser, 1996
\vskip 0.25cm
\item{[Moreno]} C. Moreno, {\it  Algebraic curves over finite fields.}
Cambridge Tracts in Mathematics, 97. Cambridge University Press,
Cambridge, 1991
\vskip 0.25cm
\item{[Neukirch]} J. Neukirch, {\it Algebraic Number Theory}, Grundlehren
der Math. Wissenschaften, Vol. {\bf 322}, Springer-Verlag, 1999
\vskip 0.25cm
\item{[Serre]} J.-P. Serre, {\it Algebraic Groups and Class Fields}, GTM
117, Springer (1988)
\vskip 0.25cm
\item{[Tate]} J. Tate, Fourier analysis in number fields and Hecke's
zeta functions, Thesis, Princeton University, 1950
\vskip 0.25cm
\item{[Weil]} A. Weil, {\it Basic Number Theory}, Springer-Verlag, 1973
\vskip 0.25cm
\item{[Weng1]} L. Weng,  Riemann-Roch Theorem, Stability and
New Zeta Functions for Number Fields,  version $1{1\over 2}$, preprint
\vskip 0.25cm
\item{[Weng2]} L. Weng,  A Program for Geometric Arithmetic, preprint
\end